\documentclass[12pt, a4paper]{amsart}

\usepackage{graphicx}
\usepackage{wrapfig}
\usepackage{comment}
\usepackage{lipsum}

\usepackage{pifont}

\usepackage[
unicode=true,
plainpages = false,
pdfpagelabels,
bookmarks=true,
bookmarksnumbered=true,
bookmarksopen=true,
breaklinks=true,
backref=false,
colorlinks=true,
linkcolor = DarkBlue,
urlcolor  = DarkBlue,
citecolor = DarkRed,
anchorcolor = green,
hyperindex = true,
hyperfigures
]{hyperref}
\hypersetup{
pdftitle={An explicit embedding of the group Q into a finitely presented group},
pdfauthor={V.H. Mikaelian},
pdfsubject={An explicit embedding of the group Q into a finitely presented group}
}

\usepackage[dvipsnames]{xcolor}
\colorlet{LightRubineRed}{RubineRed!70!}
\colorlet{Mycolor1}{green!10!orange!90!}
\definecolor{DarkRed}{HTML}{cc0000}
\definecolor{ChapterHeadColor}{HTML}{cc0000}
\definecolor{PartHeadColor}{HTML}{cc0000}
\definecolor{DarkBlue}{HTML}{0000cc}
\definecolor{QuoteColor}{HTML}{665665}

\usepackage[top=2.5cm, bottom=2.5cm, outer=2.5cm, inner=2.5cm, heightrounded
]{geometry}

\renewcommand{\S}{{\mathcal S}}
\newcommand{\Zz}{{\mathcal Z}}
\newcommand{\E}{{\mathcal E}}

\newcommand{\N}{{\mathbb N}}
\newcommand{\Z}{{\mathbb Z}}
\newcommand{\Q}{{\mathbb Q}}

\newcommand{\1}{\{1\}}

\newcommand{\bigast}{\textrm{\footnotesize\ding{91}}}

\usepackage[bitstream-charter]{mathdesign}

\DeclareSymbolFont{cmsymbols}{OMS}{cmsy}{m}{n}
\SetSymbolFont{cmsymbols}{bold}{OMS}{cmsy}{b}{n}
\DeclareSymbolFontAlphabet{\mathcal}{cmsymbols}

\theoremstyle{plain}
\newtheorem{Theorem}{Theorem}[section]
\newtheorem{Lemma}[Theorem]{Lemma}
  
\newtheorem{Corollary}[Theorem]{Corollary}

\theoremstyle{definition}

\theoremstyle{remark}
\newtheorem{Remark}[Theorem]{Remark} 
\newtheorem{Example}[Theorem]{Example} 

\numberwithin{equation}{section}

\theoremstyle{plain}
\newtheorem{innercustomgeneric}{\customgenericname}
\providecommand{\customgenericname}{}
\newcommand{\newcustomtheorem}[2]{%
\newenvironment{#1}[1]
{%
\renewcommand\customgenericname{#2}%
\renewcommand\theinnercustomgeneric{##1}%
\innercustomgeneric
}
{\endinnercustomgeneric}
}

\newcustomtheorem{TheoremABC}{Theorem}

\usepackage[bitstream-charter]{mathdesign}

\DeclareSymbolFont{cmsymbols}{OMS}{cmsy}{m}{n}
\SetSymbolFont{cmsymbols}{bold}{OMS}{cmsy}{b}{n}
\DeclareSymbolFontAlphabet{\mathcal}{cmsymbols}

\begin{document}



\title[Explicit construction of benign subgroup]{Explicit construction of benign subgroup for Higman's\\ reversing operation}

\author{V. H. Mikaelian
}

\dedicatory{To the living memory of Professor Boris Isaakovich Plotkin, whose ideas and  perseverance inspired generations}

\begin{abstract}
For the Higman  reversing operation $\rho$ and for a set of integer-valued functions $\mathcal X$ the following has been proved. 
Let the subgroup $A_{\mathcal X}$ be benign in the free group $F$, let the respective finitely presented overgroup $K_{\mathcal X}$ with its finitely generated subgroup $L_{\mathcal X}$ be given for $A_{\mathcal X}$ explicitly, and let the set $\mathcal Y = \rho\mathcal X$ be obtained from $\mathcal X$ by the reversing operation $\rho$. 
Then $A_{\mathcal Y}$ also is benign in $F$ and, moreover, the finitely presented overgroup $K_{\mathcal Y}$ with its finitely generated subgroup $L_{\mathcal Y}$ can also be given explicitly for it.  
The current work is a step in development of a general algorithm for construction of explicit Higman embeddings of recursive groups into finitely presented groups. 
\end{abstract}

\date{\today}

\maketitle

\setcounter{tocdepth}{1}

\let\oldtocsection=\tocsection
\let\oldtocsubsection=\tocsubsection
\let\oldtocsubsubsection=\tocsubsubsection
\renewcommand{\tocsection}[2]{\hspace{-12pt}\oldtocsection{#1}{#2}}
\renewcommand{\tocsubsection}[2]{\footnotesize \hspace{4pt} 
\oldtocsubsection{#1}{#2}}
\renewcommand{\tocsubsubsection}[2]{\footnotesize \hspace{34pt}\oldtocsubsubsection{#1}{#2}}

{\footnotesize \tableofcontents}

\section{Introduction}
\label{SE Introduction}

\subsection{Embeddings of recursive groups and Higman's basic operations}
\label{SU Embeddings of recursive groups and Higman's basic operations} 

This paper continues our recent research 
\cite{The Higman operations and  embeddings, 
Embeddings using universal words, 
On explicit embeddings of Q} 
on explicit embeddings of recursive groups into finitely presented groups. 
According to Higman's famous Theorem~1 in \cite{Higman Subgroups of fP groups}, a finitely generated group can be embedded into a finitely presented group if and only if it is recursive. 
According to Corollary on \textit{p.}456 in \cite{Higman Subgroups of fP groups} this theorem's requirement is sufficient for the effectively enumerable \textit{countably} generated groups: they also are embeddable into finitely presented groups, as soon as they are recursive. It is not hard to see that in this case the condition is no longer necessary though, see Example~\ref{EX Why Higman's theorem is not iff for countably generated groups} below.

Higman's result states \textit{possibility} of such an embedding without giving an \textit{explicit} algorithm for its construction. Moreover, even for the additive group of rational numbers $\Q$, which has a fairly uncomplicated recursive presentation, the question of its \textit{explicit} embedding into a finitely represented group remained open until very recently.
That problem was directly stated, in particular, in \cite{kourovka,De La Harpe 2000,Johnson on Higman's interest}, and it was solved in \cite{Belk Hyde Matucci, On explicit embeddings of Q}, 
see also our report \cite{Explicit embeddings Moscow 2018} of 2018. 

\medskip
\cite{On explicit embeddings of Q} intended to cover the embedding of the recursive group $\Q$ into some finitely presented groups, but some wide generalizations of its methods are possible. Hence, continuing the cited research we recently started generalization of arguments of \cite{On explicit embeddings of Q} to a common method that works not only for $\Q$ but also for \textit{arbitrary recursive group}, namely, for any group $G=\langle X\,|\, R \rangle$ for which: 
\begin{enumerate}
\item the set of generators $X$ is an effectively enumerable finite or countable set;
\item the set of relations $R$ is a recursively enumerable set.
\end{enumerate}
Intuitively, $G$ can be understood as a group given via $G=\langle x_1, x_2, \ldots \,|\, r_1, r_2, \ldots \rangle$, for which we possess an algorithm that outputs the relations $r_i$ of $G$ in whatever order (our algorithm recursively enumerates them). 
More precisely, by the classic definition of recursive enumeration, we assume that we are given a bijection from $R$ to a subset of the non-negative integers, and that subset can be built from the \textit{zero}, \textit{successor}, and \textit{projection} functions by means of \textit{composition}, \textit{primitive recursion}, and \textit{minimization} operations, see \cite{Davis, Rogers}.

\medskip
In order to construct the desired explicit embedding, we over the past few years, have been systematically analyzing the main stages of the Higman's construction  \cite{Higman Subgroups of fP groups}, trying to make them explicit, in case they are not already so. 

A key tool by means of which Higman concatenates the logical concept of recursion to group-theoretical concepts is the notion of benign subgroup. 
A subgroup $H$ of a finitely generated group $G$ is called a \textit{benign} subgroup in $G$ if the group $G$ can be embedded into a finitely presented group $K$ which has a finitely generated subgroup $L$ such that $G \cap L = H$ holds, see more details and references in Section~\ref{SU Benign subgroups} below. 

For a given recursive group $G$ the construction of \cite{Higman Subgroups of fP groups} first outputs a specific set $\mathcal X$ of functions from $\Z$ to $\Z$ with finite supports. 
Then Higman introduces a respective subgroup $A_{\mathcal X}$ inside the free group $F$ of rank $3$, see definition in Section~\ref{SU Defining subgroups by sets of functions} below. The group $A_{\mathcal X}$ turns out to be benign in $F$ if and only if $\mathcal X$ is recursively enumerable, see Theorem 3 and Theorem 4 in \cite{Higman Subgroups of fP groups}.

A large part of \cite{Higman Subgroups of fP groups} is dedicated to showing that if $A_{\mathcal X}$ is benign in the free group $F=\langle a,b,c\rangle$ of rank $3$, and if the set $\mathcal Y$ is obtained from $\mathcal X$ by one of nine specific \textit{Higman operations}, see \eqref{EQ Higman operations} in Section~\ref{SU The Higman operations} below, then $A_{\mathcal Y}$ also is benign in $F$. 
To build an \textit{explicit} embedding of $G$ into a finitely presented group we, in this and in following papers, for each of nine Higman operations  \eqref{EQ Higman operations} plan to \textit{explicitly} show the finitely presented group $K_{\mathcal Y}$ with its finitely generated subgroup $L_{\mathcal Y}$ admitting the equality $F \cap L_{\mathcal Y} = A_{\mathcal Y}$, as soon as $\mathcal Y$ is obtained from $\mathcal X$ by one of Higman operations. 
Also, we plan to explicitly write the overgroup $K_{\mathcal Y}$ by its generators and defining relations, and write the subgroup $L_{\mathcal Y}$ by its generators, as soon as $K_{\mathcal X}$,  $L_{\mathcal X}$ are explicitly given in that format for the set $\mathcal X$, while  $F \cap L_{\mathcal X} = A_{\mathcal X}$ holds.

\medskip
We understand that the above concise outline is not enough to give understanding of the position of our Theorem~\ref{TH Theorem for rho} in the framework of Higman's construction in \cite{Higman Subgroups of fP groups}.  
Thus, in Appendix in Chapter~\ref{SU Some of the main steps of Higman's construction} we present a more detailed depiction of the steps of \cite{Higman Subgroups of fP groups}, to spot the fragment of \cite{Higman Subgroups of fP groups}  which  Theorem~\ref{TH Theorem for rho} concerns to. Also notice Remark~\ref{RE we do not use <R>^F} and Remark~\ref{RE we do not use f a b} in which we stress some of the steps in \cite{Higman Subgroups of fP groups} which cause extra problems, and which we try avoid via modifications of the construction.

\subsection{Motivation for Higman's reversing operation $\rho$ and Theorem~\ref{TH Theorem for rho}}
\label{SU HMotivation for Higman's reversing operation rho and Theorem 1.1} 

Higman's operation $\rho$ is defined in Section 2 in \cite{Higman Subgroups of fP groups}, and it is restated with a simple example in Section~\ref{SU The Higman operations} below, see also 
Section 2.3 in  \cite{The Higman operations and  embeddings}. Notice that there is a very minor typo in \cite{Higman Subgroups of fP groups}: on \textit{p.}\,459, \textit{l.}\,10 its definition has to be not ``$i\rho=-1$'' but ``$i\rho=-i$''. 

\medskip
In the context outlined in previous section, the operation $\rho$ is interesting, we think, for the following reasons:

\textit{Firstly}, constructing the explicit embedding for the group $\Q$ we have already discussed the Higman operations 
$\sigma,\, 
\tau,\,
\zeta,\, 
\pi,\, 
\omega_m$  in \cite{On explicit embeddings of Q}. The operation $\rho$ was \textit{not} considered in \cite{On explicit embeddings of Q} because it was not needed for explicit embedding of $\Q$. 
On the other hand, the Higman operations $\iota,
\upsilon$ demand no particular discussion as the respective groups $K_{\mathcal Y}$ and $L_{\mathcal Y}$ are easily built for them by the original tools of \cite{Higman Subgroups of fP groups} already. Thus, $\rho$ is one of the few remaining operations which were \textit{not yet} properly covered so far, and for which non-evident work still is required.

\textit{Secondly}, discussion of the operation $\rho$ very clearly shows \textit{how much more} work is needed for the \textit{explicit} embedding. Check Higman's extremely brief first paragraph in the proof of Lemma 4.6 on \textit{p.}\,470 in \cite{Higman Subgroups of fP groups}. For the case $\mathcal Y=\rho \mathcal X$ that paragraph simply argues that $A_{\mathcal Y}=A_{\rho \mathcal X}$ is benign, as soon as $A_{\mathcal X}$ is benign, by showing that there is a simple isomorphism $\alpha$ sending $d_i$ to $d_{\,-i}$.
But \textit{how} can the respective groups $K_{\mathcal Y}$ and $L_{\mathcal Y}$ be written \textit{explicitly} for this case? That is done in Chapter~\ref{SU The proof for rho} below, including the explicit presentation \eqref{EQ K rho X}. 
Yes, a reader with sound knowledge in free constructions could be able to write down the group assumed by Higman, but the efforts needed for that job are far from being trivial, and $\rho$ is the first of operations \eqref{EQ Higman operations} that very well stresses stress that fact. 

\medskip
The above motivation brings us to the the main objective of this paper: prove that, if for the given set $\mathcal X$ of functions from $\Z$ to $\Z$ with finite supports, the respective subgroup $A_{\mathcal X}$ is benign in $F=\langle a,b,c\rangle$, and if the respective groups 
$K_{\mathcal X}$ and $L_{\mathcal X}$ are explicitly given for it, then $A_{\mathcal Y}=A_{\rho \mathcal X}$ also is benign in $F$, and the new groups 
$K_{\rho \mathcal X}$ and $L_{\rho \mathcal X}$ can also be written down for this updated group $A_{\mathcal Y}$:

\begin{Theorem}
\label{TH Theorem for rho}
Let $\mathcal X$ be any subset of $\E$, let $\mathcal Y = \rho\mathcal X$ be the subset obtained from $\mathcal X$ by the Higman operation $\rho$ from \eqref{EQ Higman operations},  let $F$ be the free group of rank $3$, and let $A_{\mathcal X}$ and $A_{\mathcal Y}$ be the respective subgroups in $F$.
If $A_{\mathcal X}$ is benign in  $F$, then $A_{\mathcal Y}$ is also  benign in  $F$.

Moreover, if the respective finitely presented group $K_{\mathcal X}$ with its finitely generated subgroup $L_{\mathcal X}$
are explicitly given for $A_{\mathcal X}$, then  
$K_{\mathcal Y}$ and $L_{\mathcal Y}$ can also be explicitly given for $A_{\mathcal Y}$.
\end{Theorem}

Figure~\ref{FI Passage from the benign subgroup} below illustrates the passage from $\mathcal X$ and $A_{\mathcal X}$ to $\mathcal Y=\rho \mathcal X$ and to $A_{\mathcal Y} = A_{\rho \mathcal X}$ in this 
theorem. Its proof will be given in Chapter~\ref{SU The proof for rho} below, after some definitions and preparations in  
Chapter~\ref{SE Preliminary notation and references} and Chapter~\ref{SE Additional constructions for the proof of Theorem}.

\begin{figure}[h]
\includegraphics[width=300px]{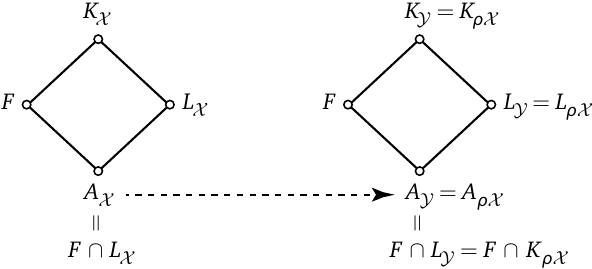}
\caption{Passage from $A_{\mathcal X}$ to the new benign subgroup $A_{\mathcal Y} = A_{\rho \mathcal X}$ in Theorem~\ref{TH Theorem for rho}.} 
\label{FI Passage from the benign subgroup}
\end{figure}

We would like to announce that the work on analogs of Theorem~\ref{TH Theorem for rho} for the remaining few types of Higman operation \eqref{EQ Higman operations} currently is in progress, and hence construction of the explicit Higman embedding algorithm seems to be within reach.   Notice, however, that such an explicit embedding is \textit{not} automatically following for the promised analogs of Theorem~\ref{TH Theorem for rho} for all operations alone, because Higman's construction contains other steps to take care of, besides those in  sections 3 and 4 in \cite{Higman Subgroups of fP groups}. 
See Chapter~\ref{SU Some of the main steps of Higman's construction} at the end of the current paper, with some more particularized description of Higman's steps in \cite{Higman Subgroups of fP groups}. 

\medskip 
Since we wish to have as few as possible repetitions of the fragments of proof parts from our previous work, we list in 
Chapter~\ref{SE Preliminary notation and references}
and
Chapter~\ref{SE Additional constructions for the proof of Theorem} 
the main tools and auxiliary facts along with references to details in~\cite{
Embeddings using universal words,
The Higman operations and  embeddings, 
Auxiliary free constructions for explicit embeddings, 
On explicit embeddings of Q,
A modified proof for Higman}. 
%

%

\subsection*{Acknowledgements}
\label{SU Acknowledgements}
I am very much thankful to the Referee for highly valuable remarks. In particular, his suggestion to add Chapter~\ref{SU Some of the main steps of Higman's construction} helps to stress the motivation of this article, and to reason why we are considering Higman's  reversing operation $\rho$.

The current work is supported by the grant 25RG-1A187 of SCS MES Armenia.

I am thankful to the Fachliteratur Program of the German Academic Exchange Service DAAD for the academic literature provided over the past years, grant A/97/13683.

\section{Preliminary notation and references}
\label{SE Preliminary notation and references} 

\subsection{Integer-valued functions $f$}
\label{SU Integer functions f} 

Following \cite{Higman Subgroups of fP groups} denote by $\mathcal E$ the set of \textit{all} functions $f : \Z \to \Z$ with finite supports; 
as usual the support of $f$ is the set ${\rm sup}(f)=\big\{i\in \Z \;|\; f(i)\neq 0\big\}$.
%
When $f$ for a certain fixed $m=1,2,\ldots$ has the property that 
$f(i)=0$ for \textit{all} $i<0$ and $i\ge m$, then it is comfortable to record it as a sequence $f=(j_0,\ldots,j_{m-1})$, where $f(i)=j_i$ for $i=0,\ldots,m-1$. For example, the function $f$ which sends the integers $0, 1, 2$ to $2, 5, 3$,\; and all other integers to $0$, is recorded as $f=(2, 5, 3)$.
Clearly, $m$ may not be uniquely defined for $f$\!, and where necessary we may add extra zeros at the end of a sequence, e.g., the previous function can be also recorded as $f=(2, 5, 3,0,0)$. 
%
See more on such functions in  Section 2.2 in \cite{The Higman operations and  embeddings}.
Denote $\E_m\!=\big\{
(j_0,\ldots,j_{m-1})
\;\mathrel{|}\;
j_i\!\in \!\Z,\; i\!=\!0,\ldots,m-1\big\}$ for $m=1,2,\ldots$ 
Clearly $\E_m \subseteq \E$, and by the agreement above $\E_{m'}\subseteq \E_m$,
in case $m' \le m$.

For any $f\in \mathcal E$ and $j\in \Z$  we define the function $f_{j}^+$\! as follows: 
$f_{j}^+(i)=f(i)$ for all $i\!\neq\! j$, and   
$f_{j}^+(j) = f(j)+\!1$.
When $f\in \E_m$, we define $f^+ = f_{m-1}^+$, say, for the above $f\!=\!(2, 5, 3)\in\E_3$ we have $f_{1}^+\!\!=(2, 6, 3)$ and 
$f^+\!\!=(2, 5, 4)$.

\subsection{The Higman operations}
\label{SU The Higman operations} 
Higman defines the following specific  operations that  transform the subsets of $\E$ to some new subsets of  $\E$:
\begin{equation}
\tag{H}
\label{EQ Higman operations}
\iota,\; 
\upsilon,\; 
\rho,\; 
\sigma,\; 
\tau,\; 
\theta,\; 
\zeta,\; 
\pi,\; 
\omega_m,
\end{equation}
$m=1,2,\ldots$\,,
see Section 2 in  \cite{Higman Subgroups of fP groups} or Section 2.3 in  \cite{The Higman operations and  embeddings}  for details.

In particular, one of the Higman operations \eqref{EQ Higman operations} is the \textit{reversing} operation $\rho$ defined as follows: 
for any subset $\mathcal X \subseteq \mathcal E$ the set 
$\mathcal Y = \rho\mathcal X$ consists of all functions $f\in \mathcal E$ for which there is a $g\in \mathcal X$ such that $f(i)=g(-i)$ for all $i\in \Z$. Intuition behind this definition is uncomplicated: $\rho$ just reverses the ``coordinates'' of $f$ symmetrically about the origin $0$.
Say, if $f$ is the function sending 
the integers $-1, 0, 1, 2$ to $3,2,9,8$,\; and all other integers to $0$,
then $\rho f$ is sending $-2, -1, 0, 1$ to $8, 9, 2, 3$,\; and all other integers to $0$.

\subsection{Defining subgroups by sets of functions}
\label{SU Defining subgroups by sets of functions}

Let $F=\langle a, b, c\rangle$ be a free group of rank $3$. For any $i\in \Z$ denote $b_i=b^{c^i} = c^{-i} b \, c^i$. Then for any function $f\! \in \E$ one can define in $F$ the following products: 
\begin{equation}
\label{EQ defining b_f and a_b}
b_{f} = \cdots
b_{-1}^{f(-1)}
b_{0}^{f(0)} 
b_{1}^{f(1)}\cdots 
\quad \text{and} \quad
a_{f} = a^{b_{f}} = b^{-1}_{f} a\, b_{f}.
\end{equation}
Say, for the above function 
$f$ sending\,  $-1, 0, 1, 2$ to $3,2,9,8$ we have 
$b_f\!=\!b_{-1}^{3}
b_{0}^{2} 
b_{1}^{9}
b_{2}^{8}$.
When $f$ is in $\E_m$, we can more comfortably record it as 
$f=(j_0,\ldots,j_{m-1})$, see Section~\ref{SU Integer functions f}, and then write
$b_{f} = b_0^{j_0}  \cdots b_{m-1}^{j_{m-1}}$.
Say, for $f\!=\!(2, 5, 3)$
we may put
$b_f\!=\!
b_{(2, 5, 3)}\!=
b_{0}^{2} 
b_{1}^{5}
b_{2}^{3}$ and
$a_{(2, 5, 3)}\!=a^{b_{(2, 5, 3)}}$.

For a set $\mathcal X$ of functions from $\E$ denote by $A_{\mathcal X}$ the subgroup generated in $F$ by all the conjugates $a_{f} = a^{b_{f}}$ with $f\!\in \mathcal X$.
Say, $A_{\E_2}$ is the subset generated in $F$ by all words of type
$a^{b_{(j_0,\; j_1)}}
=  
b_{1}^{\!-j_1}
b_{0}^{\!-j_0}
\cdot a \cdot b_{0}^{j_0} 
b_{1}^{j_1}$ for all possible integers $j_0, j_1 \in \Z$.

Applying the operation $\rho$ to the above function $f$, we obtain from the elements
$b_f = b_{\!-1}^3 b_{0}^2 b_{1}^ 9 b_{2}^8$,\;
$a_f=a^{b_{\!-1}^3 b_{0}^2 b_{1}^ 9 b_{2}^8}$ \; the following elements:\;
$b_{\rho f} = b_{\!-2}^8 b_{\!-1}^9 b_{0}^ 2 b_{1}^3$,\;
$a_{\rho f}=a^{b_{\!-2}^8 b_{\!-1}^9 b_{0}^ 2 b_{1}^3}$.

For some technical purposes we may use the above notation with \textit{other} free generators also. Say, in the free group $\langle d,e\rangle$
we may set  
$d_i=d^{e^i}$\!,\;\;
$d_{f} = \cdots
d_{-1}^{f(-1)}
d_{0}^{f(0)} 
d_{1}^{f(1)}\cdots $\;
Or, we may take an isomorphic copy 
$\bar F=\langle \bar a, \bar b, \bar c\rangle$ of $F=\langle a, b, c\rangle$, and then use in $\bar F$ the elements $\bar b_i, \bar b_f, \bar a_f$ as expected.

\subsection{Benign subgroups}
\label{SU Benign subgroups} 

For detailed information on benign subgroups, quickly mentioned in the Introduction, we refer to  Sections 3, 4 in \cite{Higman Subgroups of fP groups}, see also Section 3 in \cite{A modified proof for Higman}.
If the group $G$ for its benign subgroup $H$ can be embedded into a finitely presented group $K$ with a finitely generated subgroup $L$ such that $G \cap L = H$,
then we may stress this correlation by denoting   $K=K_H$ and $L=L_H$.
If $A_{\mathcal X}$ is benign in $F$, we may denote the respective groups via
$K_{\mathcal X}$
and
$L_{\mathcal X}$ to have comfortable notation.

\begin{Remark}
\label{RE finite generated is benign}
It is very easy to see that arbitrary \textit{finitely generated} subgroup $H$ of any \textit{finitely presented} group $G$ is benign in $G$.  
Indeed, 
$H = K_H \cap L_H$ for the choice $K_H = G$ with $L_H = H$.
We are going to often use this remark in the sequel.
\end{Remark}

\subsection{Free constructions}
\label{SU Free constructions} 

For background information on free products with amalgamation and on HNN-extensions we refer to 
\cite{Bogopolski} and \cite{Lyndon Schupp}. 
Since notation varies  in the literature, 
we fix them for this paper as follows.
Suppose that $G$ and $H$ are two groups, and $A$ and $B$ are isomorphic subgroups of $G$ and $H$, respectively, with an isomorphism $\varphi : A \to B$. Then the generalized free product of $G$ and $H$ with amalgamation of $A$ and $B$ by $\varphi$ is denoted by $G*_{\varphi} H$.
When $G$ and $H$ are overgroups of the same subgroup $A$, and $\varphi: A \to A$ is the isomorphism fixing any element $A$, then we denote the corresponding amalgamated product as $G*_{A} H$.

If $G$ has subgroups $A$ and $B$ isomorphic under  $\varphi : A \to B$, then the corresponding HNN-extension  of the base group $G$ 
by some stable letter $t$
with respect to 
$\varphi$ is denoted by
$G*_{\varphi} t$.
In case when $A=B$ and $\varphi$ is the isomorphism just fixing $A$, we denote the corresponding HNN-extension by $G*_{A} t$.
We also use HNN-extensions $G *_{\varphi_1, \varphi_2, \ldots} (t_1, t_2, \ldots)$ and 
$G *_{A_1,\; A_2, \ldots} (t_1, t_2, \ldots)$  
with more than one stable letters, see \cite{A modified proof for Higman} for details.

Below we are going to use a series of facts about certain specific subgroups in  free constructions $G*_{\varphi} H$, \, $G*_{A} H$, \,$G*_{\varphi} t$ and $G*_{A} t$. We have collected them in Section 3 of  \cite{Auxiliary free constructions for explicit embeddings} to refer to that section whenever needed.

\section{Additional constructions for the proof of Theorem~\ref{TH Theorem for rho}
}
\label{SE Additional constructions for the proof of Theorem}

\noindent
In this section we amass some constructions, definitions, results, and links from \cite{Embeddings using universal words, The Higman operations and  embeddings, Auxiliary free constructions for explicit embeddings, On explicit embeddings of Q} to refer to them wherever needed.  
The main aim here is to recall the definition of the group $\mathscr{A}$ from \cite{On explicit embeddings of Q}. That group will play a crucial role in this paper. Therefore, we give an explicit finite presentation of $\mathscr{A}$, and an example that will give some experience for certain calculations in this group. However, we believe that some general constructions that lead to the group $\mathscr{A}$ are important and, therefore, we explain them also before turning to the construction of $\mathscr{A}$.

\subsection{The $\bigast$-construction}
\label{SU The *-construction} 

Let $G, M, K_1,\ldots,K_r$ be a family of groups meeting the following conditions: 
\begin{equation}
\label{EQ conditions for star construciton} 
\text{$G\le M  \le K_1,\ldots,K_r$, \;and\; $K_i \cap\,  K_j=M$\; for any distinct\; $i,j=1,\ldots,r$.}
\end{equation} 
Taking an arbitrary subgroup $L_i$ in $K_i$ for each $i=1,\ldots,r$,
we in Section~4 of
\cite{Auxiliary free constructions for explicit embeddings} built an auxiliary technical ``$\bigast$-construction'' applying  a nested combination of HNN-extensions and amalgamations:
\begin{equation}
\label{EQ nested Theta  multi-dimensional short form} 
\textstyle
\bigast_{i=1}^{r}(K_i, L_i, t_i)_M\!
= \Big(\! \cdots
\Big( \big( (K_1 *_{L_1} t_1) *_M\! (K_2 *_{L_2} t_2) \big) *_M\!  (K_3 *_{L_3} t_3)\Big)\cdots 
\!\Big) *_M\!  (K_r *_{L_r} t_r).
\end{equation} 
Recall that $K_i *_{L_i} t_i$ is the HNN-extension of the group $K_i$ with the base subgroup $L_i$ and with the stable letter $t_i$ fixing $L_i$ elementwise, $i=1,\ldots,r$.
Then \eqref{EQ nested Theta  multi-dimensional short form} is the free product
of these HNN-extensions $K_i *_{L_i} t_i$, all  amalgamated in their \textit{common} subgroup $M$.

\vskip5mm

\begin{figure}[h]
\includegraphics[width=300px]{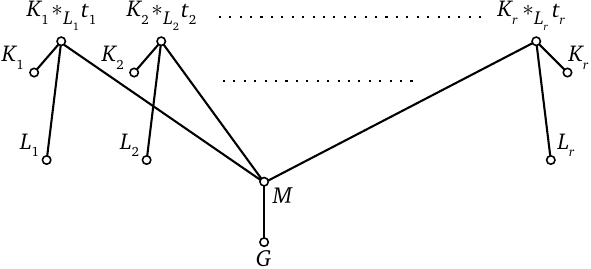}
\caption{Construction of the group\, $\bigast_{i=1}^{r}(K_i, L_i, t_i)_M$.} 
\label{FI Construction of the group *}
\end{figure}

\medskip
In some extreme cases the $\bigast$-construction \eqref{EQ nested Theta  multi-dimensional short form}, in fact, coincides with certain well known ``classic'' constructions.

\begin{Example}
When a certain group $G$ and its subgroups $A_1,\ldots,A_r$ are fixed, then taking $K_i=M=G$ and $L_i=A_i$ for each $i=1,\ldots,r$,  we 
have the HNN-extensions $G *_{A_i} t_i$, and then 
$\bigast_{i=1}^{r}(G, A_i, t_i)_G$ is the amalgamation of all such $G *_{A_i} t_i$ by their subgroup $G$. 
In this case the $\bigast$-construction is the usual HNN-extension with multiple stable letters $t_1,\ldots,t_r$: 
$$
G *_{A_1,\ldots,\,A_r } \!(t_1,\ldots,t_r). 
$$
\end{Example} 

\begin{Example}
If we again take a random group $G$, and put $K_i=M=G$ for each $i=1,\ldots,r$, and choose \textit{trivial} subgroups  $L_i=A_i=\1$, then $\bigast_{i=1}^{r}(G, A_i, t_i)_G$ simply is the ordinary free product of $G$ with the free group $\langle t_1,\ldots,t_r \rangle\cong F_r$ of rank $r$:
$$
G * \langle t_1,\ldots,t_r \rangle 
= G * F_r.
$$ 
\end{Example} 

\begin{Example}
Consider the case when $L_i=K_i=M=G$ for each $i=1,\ldots,r$.  Then in $G *_{G} t_i$ conjugation by stable letter $t_i$ just fixes the whole $G$, which means this HNN-extension is the direct product $G \times \langle t_i \rangle$ for $\langle t_i \rangle \cong \Z$. Then the $\bigast$-construction $\bigast_{i=1}^{r}(G, G, t_i)_G$ turns out to be  the \textit{direct} product: 
$$
G \times \langle t_1,\ldots,t_r \rangle
\cong
G \times F_r.
$$ 
Factorizing the above with the relations $[t_i, t_j]$, with $i,j=1,\ldots,r$, we get an even simpler direct product:
$$
G \times \langle t_1,\ldots,t_r \rangle
\cong
G \times \Z^r.
$$ 
\end{Example}

\medskip
The following lemmas \ref{LE intersection in bigger group multi-dimensional}\,--\,\ref{LE join in HNN extension multi-dimensional} are proven in sections 3 and 4 in \cite{Auxiliary free constructions for explicit embeddings}, and in each of them the groups $G,\; M,\; K_i,\; L_i,\; A_i$ are those mentioned in the definition of the $\bigast$-construction \eqref{EQ nested Theta  multi-dimensional short form}. First, we discover a familiar HNN-extension of $G$ inside \eqref{EQ nested Theta  multi-dimensional short form}:

\begin{Lemma}
\label{LE intersection in bigger group multi-dimensional}
In the above notation the following equality holds in ${\bigast}_{i=1}^{r}(K_i, L_i, t_i)_M$:
$$
\langle G, t_1,\ldots,t_r \rangle= 
G *_{A_1,\ldots,\,A_r } \!(t_1,\ldots,t_r).
$$
\end{Lemma}

\subsection{Intersections and joins of benign subgroups}
\label{SU Intersections and joins of benign subgroups} 

It turns out that the intersection $I$ of the subgroups $A_i$ from the previous section can be obtained as intersection of \textit{just two} conjugate subgroups in \eqref{EQ nested Theta  multi-dimensional short form}:

\begin{Lemma}
\label{LE intersection in HNN extension multi-dimensional}
In the above notation, if $I=\bigcap_{\,i=1}^{\,r} \,A_i$, then in $G *_{A_1,\ldots,\,A_r} (t_1,\ldots,t_r)$ and, hence, also in ${\bigast}_{i=1}^{r}(K_i, L_i, t_i)_M$, we have
\begin{equation}
\label{EQ gemeral intersection in HNN}
\textstyle
G \cap G^{t_1 \cdots\, t_r}
= I.
\end{equation}
\end{Lemma}

The approximate analog of the above lemma holds for the join $J$ of the subgroups $A_i$:

\begin{Lemma}
\label{LE join in HNN extension multi-dimensional}
In the above notation, if 
$J=\big\langle\bigcup_{\,i=1}^{\,r} \,A_i\big\rangle$, then in $G *_{A_1,\ldots,\,A_r} (t_1,\ldots,t_r)$ and, hence, also in ${\bigast}_{i=1}^{r}(K_i, L_i, t_i)_M$, we have
\begin{equation}
\label{EQ gemeral intersection in HNN multi-dimensional}
\textstyle 
G \cap \big\langle 
\bigcup_{\,i=1}^{\,r} \,G^{t_i} \big\rangle
=J.
\end{equation}
\end{Lemma}


\smallskip
These lemmas allow us to build new benign subgroups from the existing ones:

\begin{Corollary}
\label{CO intersection and join are benign multi-dimensional}
If the subgroups $A_1,\ldots,\,A_r$ are benign in a finitely generated group $G$, then
\begin{enumerate}
\item 
\label{PO 1 CO intersection and join are benign multi-dimensional}
their intersection $I=\bigcap_{\,i=1}^{\,r} \,A_i$ also is benign in $G$;
\item 
\label{PO 2 CO intersection and join are benign multi-dimensional}
their join $J=\big\langle\bigcup_{\,i=1}^{\,r} \,A_i\big\rangle$ also is benign in $G$.
\end{enumerate}
Moreover, if the finitely presented groups $K_i$ with their finitely generated subgroups $L_i$ can be 
given for each $A_i$ explicitly, then the respective finitely presented overgroups $K_I$ and $K_J$ with finitely generated  subgroups $L_I$ and $L_J$
can also be given for $I$ and for $J$ explicitly.
\end{Corollary}

As finitely presented overgroups $K_I$ and $K_J$ promised in this corollary one may take the above mentioned $\bigast$-construction $\textstyle{\bigast}_{i=1}^{r}(K_i, L_i, t_i)_M$ which evidently is \textit{finitely presented}, as soon as $K_1,\ldots,K_r$ are finitely presented, and the subgroups $M,\; L_1,\ldots,L_r$ are finitely generated.
Also, one may choose the subgroups $L_I\!=G^{t_1 \cdots\, t_r}$ and $L_J=\big\langle 
\bigcup_{\,i=1}^{\,r} \,G^{t_i} \big\rangle$ which, evidently, are \textit{finitely generated}. 
See Figure~6 in \cite{Auxiliary free constructions for explicit embeddings} illustrating this corollary.

\subsection{The group $\mathscr{A}$ from \cite{On explicit embeddings of Q}}
\label{SU The group A}

Let $\langle b,c \rangle$ be a free group of rank $2$ and let $m$ be an arbitrary integer. In \cite{Auxiliary free constructions for explicit embeddings} we defined two injective homomorphisms $\xi_m$ and $\xi'_m$ from $\langle b,c \rangle$ to itself by the following rules:
$$
\text{$\xi_m(b)=b_{-m+1},\;\;\;
\xi'_m(b)=b_{-m}$,\;\;\; 
$\xi_m(c)=\xi'_m(c)=c^2$\!,}
$$
see Section~\ref{SU Defining subgroups by sets of functions} above for the definition of $b_i$.
Using them we introduced in \cite{Auxiliary free constructions for explicit embeddings} the following HNN-extension to use it in various contexts: 
\begin{equation}
\label{EQ Ksi two defined}
\Xi_m = \langle b,c \rangle *_{\xi_m, \xi'_m} (t_m, t'_m).
\end{equation}
In the current paper we are going to use $\Xi_{\,1}$ for the value $m=1$ only. Thus, we have:
$$
\Xi_1 =
\big\langle b, c, t_1, t'_1
\mathrel{\;|\;} 
b^{t_1}\!=\! b,\;
b^{t'_1}\!=\! b^{c^{-1}}\!\!\!\!,\;
c^{t_1}\!=\!c^{t'_1}\!=\!c^2
\big\rangle.
$$
In \cite{Auxiliary free constructions for explicit embeddings} we have seen that the countably generated subgroup
$\langle b_1, b_2,\ldots \rangle$ is benign in $F=\langle a,b,c \rangle$ for the finitely presented overgroup $\langle a \rangle * \Xi_1$, and for its finitely generated subgroup $\langle b_1, t_1, t'_1\rangle$.
Also the subgroup $\langle a, b_{0}, b_{-1},\ldots\rangle$ is benign in $F$ for the same finitely presented group $\langle a \rangle * \Xi_1$ and for its finitely ge\-nerated subgroup $\langle a, b_{0}, t_1, t'_1\rangle$.
We use them to build the finitely presented
$\bigast$-construction:
\begin{equation}
\label{EQ Defining C}
\mathscr{C} 
= \Big( \big(\langle a \rangle\! * \Xi_1\big) *_{\langle b_1, t_1, t'_1 \rangle} u_1 \Big)
\,*_{\langle a \rangle *\, \Xi_1}
\Big(\big(\langle a \rangle\! * \Xi_1\big) *_{\langle a, b_{0}, t_1, t'_1 \rangle} u_2 \Big).
\end{equation}

The subgroup of $\mathscr{C}$, generated by $F^{u_1}$ and $F^{u_2}$ is isomorphic to the free product $F^{u_1} *\, F^{u_2}$, see Section 4.4 in \cite{
Auxiliary free constructions for explicit embeddings}. We will identify this subgroup with that free product. Let
$$
\omega: F^{u_1} *\, F^{u_2} \to F^{u_1} *\, F^{u_2}
$$
be the injective homomorphism which is the identity map on $F^{u_1}$, and is the conjugation by $b^{u_2}$ on $F^{u_2}$. Let
$\delta:F \to F$ be the injective homomorphism of $F$ sending $a, b, c$ to $a, b^c\!,\, c$.
With this notation, we defined in \cite{On explicit embeddings of Q} the following group $\mathscr{A}$ that will play a key technical role in the current paper:
$$
\mathscr{A} = \mathscr{C} \! *_{\omega, \delta}\!(d,e).
$$
This group $\mathscr{A}$ can be explicitly given by the following finite presentation:
\begin{equation}
\label{EQ relations A}
\begin{split}
\mathscr{A} 
& \!=\! 
\big\langle a, b, c, t_1, t'_1, u_1, u_2, d,e \mathrel{\;|\;}  
b^{t_1}\!=\! b,\;
b^{t'_1}\!=\! b^{c^{-1}}\!\!\!\!,\;
c^{t_1}\!=\!c^{t'_1}\!=\!c^2; \\[-2pt]
& \hskip11mm 
\text{$u_1$ fixes $b^c\!, t_1, t'_1$};\;\;\;\;
\text{$u_2$ fixes $a, b, t_1, t'_1$}; \\[-2pt]
& \hskip11mm 
\text{$d$ fixes 
$\{a,b,c\}^{ u_1}$};  \\[-2pt]
& \hskip11mm \text{$d$ sends $a^{u_2}\!\!,\; b^{u_2}\!\!,\; c^{u_2}$ to $a^{b u_2}\!,\; b^{u_2}\!,\; c^{b u_2}$};\\[-2pt]
& \hskip11mm 
\text{$e$\; sends $a,b,c$ \;to\; $a,b^c\!,\; c$}
\big\rangle,
\end{split}
\end{equation}
with ``$d$ fixes 
$\{a,b,c\}^{ u_1}$'' meaning 
$a^{ u_1 d}\!=a^{ u_1}$,\,
$b^{ u_1 d}\!=b^{ u_1}$,\,
$c^{ u_1 d}\!=c^{ u_1}$,
and with ``$d$ sends $a^{u_2}$ to $a^{b u_2}$''
meaning 
$a^{u_2 d}=a^{b u_2}$.
For later purposes we denote the set of generators of $\mathscr{A}$ by 
\begin{equation}
\label{EQ generators of XA}
X_{\! \mathscr{A}}=\big\{
a, b, c, t_1, t'_1, u_1, u_2, d,e
\big\},
\end{equation}
and the set of defining relations of $\mathscr{A}$ from \eqref{EQ relations A} by 
$R_{\! \mathscr{A}}$.
Then the presentation 
$\mathscr{A}=\langle
\,X_{\! \mathscr{A}}
\;|\;
R_{\! \mathscr{A}}
\rangle$ has $9$ generators and $2+2+3+4+3+3+3=20$ relations.

\subsection{Computing the conjugation of $a_f$ by $d_j$ in $\mathscr{A}$}
\label{SU Computing the conjugation of by d j} 

A useful computational feature takes place in $\mathscr{A}$. 
Namely, using the earlier notation 
$f_{j}^+$\!,\, 
$f_{j}^-$\!,\, 
$f^+$\!,\, 
$f^-$ 
from Section~\ref{SU Integer functions f}, 
we can for a given $a_f$ also use the elements, say, $a_{f_{j}^+}$ or $a_{f^-}$ in $\mathscr{A}$, \textit{inside} $F$.
Also, using the remark about $\langle d,e\rangle$ in Section~\ref{SU Defining subgroups by sets of functions} we can use the elements $d_i$ and $d_f$ in $\mathscr{A}$, \textit{outside} $F$.

The following lemma says that the elements $a_f$ and $a_{f_{j}^+}$ and $a_{f_{j}^-}$ are conjugated in $\mathscr{A}$ (not in $F$) by certain elements of $\mathscr{A} \backslash F$ that depend on $j$ and do not depend on $f$. This was the main reason for the construction of $\mathscr{A}$.

\begin{Lemma}
\label{LE action of d_m on f} 
For any $f \in \mathcal E$ and any $j\in \Z$ 
we have 
$
a_f^{d_j} =\! a_{f_{j}^+}$ 
and
$
a_f^{\,d_j^{-1}}\!\! = a_{f_{j}^-}
$ in $\mathscr{A}$.
\end{Lemma}

The proof to this lemma is given in \cite{A modified proof for Higman, On explicit embeddings of Q}. To gain some experience in calculations within $\mathscr{A}$, we give a simple example showing how $d_j$ acts by conjugation on $a_f$ for some particular $f$ and $j$.

\begin{Example}
\label{EX conjugation action of d_j onto a_f}
Take the sequence, say, $f=(2,5,3)$ with the product $b_f=b_0^{2}\,b_1^5\,b_2^{3}$, fix the value $j=1$, and then try to find the conjugates 
$\big(b_0^{2}\big)^{d_1}$\!\!,\;
$\big(b_1^{5}\big)^{d_1}$\! and
$\big(b_2^{3}\big)^{d_1}$ for $d_j=d_1$. 

Inside the group $F$ the isomorphism $\omega$ (and, hence, conjugation by the element $d$) leaves all the elements $ b_1, b_2,\ldots $ fixed, and sends the elements $a, b_{0}, b_{-1},\ldots$ to their conjugates $a^b\!,\, b_{0}^b,\, b_{-1}^b,\ldots$\, Hence, also taking into account the action of $\delta$ (that is, the conjugation by the element $e$) we have:
\begin{equation*}
\begin{split}
b_0^{d_1} & 
= \big(c^{-0} b c^0\big)^{e^{-1} d e}\!\!
= \big(b^{e^{-1}}\big)^{d e}\!\!
= \big(c b c^{-1}\big)^{d e}\!\!
= \big(b_{-1}\big)^{d e}\!\!
= \big(b^b_{-1}\big)^{e}\\
& 
= \big(b^{-1} c b c^{-1} b\big)^{e}
= c^{-1}b^{-1} c \cdot c \cdot c^{-1} b c \cdot c^{-1} \!\cdot c^{-1} b c
= b_1^{-1}\! \cdot b_0 \cdot b_1,
\end{split}
\end{equation*}
from where, clearly, $\big(b_0^{2}\big)^{d_1}\!=b_1^{-1} b^2_0 b_1$. 
Similarly:  
\begin{equation*}
\begin{split}
b_1^{d_1} & 
=\big(c^{-1} b c^1\big)^{e^{-1} d e}\!\!
= \big(c^{-1} b^{e^{-1}} c\big)^{d e}\!\!
= \big(c^{-1} c b c^{-1} c\big)^{d e}\!\!
= \big(b_{0}\big)^{d e}\!\!
= \big(b^b_{0}\big)^{e}\\
& 
= b^{e} 
= b^c=b_1 
= b_1^{-1}  \cdot b_1 \cdot b_1,
\end{split}
\end{equation*}
from where, $\big(b_1^{5}\big)^{d_1}\!=b_1^{-1} b^5_1 b_1$. 
The situation is a little different for the \textit{last} factor $b_2^{3}$:
\begin{equation*}
\begin{split}
b_2^{d_1} & 
=\big(c^{-2} b c^2\big)^{e^{-1} d e}\!\!
= \big(c^{-2} b^{e^{-1}} c^2 \big)^{d e}\!\!
= \big(c^{-2} c b c^{-1} c^2\big)^{d e}\!\!
= \big(b_{1}\big)^{d e}\!\!
= \big(b_{1}\big)^{e}\\
& 
= \big(c^{-1} b c\big)^{e}
= c^{-1} b^c c = b_2,
\end{split}
\end{equation*}
from where, again, $\big(b_2^{3}\big)^{d_1}\!= b^3_2$. 
Finally, we have
$$
a^{d_1} 
= a^{e^{-1} d e}\!
= \big(a^{e^{-1}}\big)^{de}\!
= \big(a^{d}\big)^{e}\!
= \big(b^{-1} a b \big)^{e}\!
= b^{-c}  a b^c
= b_1^{-1}  a b_1.
$$
Substituting the above calculated values into $b_f$, and then into $a_f$ for $f=(2,5,3)$, we get:
\begin{equation*}
\begin{split}
a_f^{d_1} & =
\big(
b_2^{\!-3}b_1^{\!-5}b_0^{\!-2}
\; a\;
b_0^{2}b_1^{5}b_2^{3} \,
\big)^{d_1}
\!\! \\
& =
\big(b_2^{\!-3} \big)^{d_1} 
\big(b_1^{\!-5} \big)^{d_1}
\big(b_0^{\!-2} \big)^{d_1}
\; (a)^{d_1} \;
\big(b_0^{2}\big)^{d_1}
\big(b_1^{5}\big)^{d_1}
\big(b_2^{3} 
\big)^{d_1}
\!\! \\
& = 
b_2^{-3}\,
\big(b_1^{\!-1} b_1^{\!-5} b_1^{\vphantom8}\big)
\big(b_1^{\!-1} b_0^{\!-2} b_1^{\vphantom8}\big)
\,\, \big(b_1^{\!-1} a b_1^{\vphantom8}\big) \,
\big(b_1^{\!-1} b_0^{2}b_1^{\vphantom8}\big)\,
\big(b_1^{\!-1} b_1^{5}b_1^{\vphantom8}\big)
\,b_2^{3}\\
& = b_2^{\!-3}
\big(b_1^{\!-1} b_1^{\!-5}\big)
b_0^{\!-2}
\;a\;
b_0^{2}
\big( b_1^{5} b_1^{1}\big)
b_2^{3}
\;=\; b_2^{\!-3}b_1^{\!-\,6}  b_0^{\!-2}
\;a\;
b_0^{2}b_1^6b_2^{3}\\
&= a_{f_{1}^+}
\end{split}
\end{equation*}
for the sequence
$f_{1}^+ \!= (2,\,\boldsymbol{5\!+\!1}\,,3)= (2,\boldsymbol{6},3)$.
Next, taking, say, $j=2$  we would have: 
$$
a_f^{d_2}=a_{f_2^+}\!=a_{f^+}
$$ 
where 
$f_{2}^+ \!= f^+ \!\!= (2,5,\,\boldsymbol{3\!+\!1})= (2,5,\boldsymbol{4})$.

\smallskip 
Hopefully, the calculation routine in the displayed example does not entomb the simple meaning of Lemma~\ref{LE action of d_m on f}: the conjugation by $d^j$ just ``lifts''  by $1$ the exponent of the factor corresponding to the $j$'th coordinate of $f$ inside $a_f$.
\end{Example}

\begin{Remark}
\label{RE order of d_i does not matter}
The following feature of this lemma will be used repeatedly. 
The \textit{order} of elements $d_i$ acting on $a_f$ does \textit{not} matter, i.e., $a_f^{d_{j_1}d_{j_2}}$\! and $a_f^{d_{j_2}\,d_{j_1}}$\! are equal for any $f, j_1, j_2$. Say, for the above $f=(2,5,3)$ we have
$$a_f=a_{(2,5,3)}^{d_{1}d_{2}}
= a_{(2,5,3)}^{d_{2}d_{1}}
=a_{(2,\boldsymbol 6, \boldsymbol 4)}.$$
\end{Remark}

\subsection{An auxiliary copy $\bar{\mathscr{A}}$ of $\mathscr{A}$}
\label{NOR SU Auxiliary copy of A built here}   

In analogy to the generating set $X_{\!\mathscr{A}}$ defined in \eqref{EQ generators of XA} we introduce a new set:  
\begin{equation}
\label{EQ generators of X bar A}
X_{\! \bar {\mathscr{A}}}=\big\{
\bar a,  \bar b,  \bar c,  \bar t_1,  \bar t'_1, \bar  u_1,  \bar u_2,  \bar d, \bar e
\big\},
\end{equation}
and using it we construct a copy $\bar {\mathscr{A}}$ of the group $\mathscr{A}$ applying the same procedure as in Section~\ref{SU The group A}.
This group has the relations $R_{\! \bar {\mathscr{A}}}$ obtained from the relations $R_{\! {\mathscr{A}}}$ of \eqref{EQ relations A} by just adding bars on each letter, such as 
$\bar b^{\bar t_1}\!=\! \bar b$,  
$\bar b^{\bar t'_1}\!=\! \bar b^{\bar c^{\; -1}}$\!\!\!,\,  etc.
In particular, inside $\bar {\mathscr{A}}$ the sub\-group
$\bar F \!=\! \langle \bar a, \bar c, \bar c\rangle$ is a free group of rank $3$. 

Next, in addition to the elements $b_i, b_f, a_f \!\in\! F$,\,
$d_i, d_f \!\in \!{\mathscr{A}}$ we  introduce the elements   
$\bar b_i, \bar b_f, \bar a_f\! \in\! \bar F$, \,
$\bar d_i, \bar d_f \!\in\! \bar {\mathscr{A}}$ 
expectedly defined as $\bar b_i\! = \bar b^{\bar c^i}$\!\!,\;
$\bar b_f\! 
= \cdots
\bar b_{-1}^{f(-1)}
\bar b_{0}^{f(0)} 
\bar b_{1}^{f(1)}\cdots$,
\,
$\bar a_f\!=\bar a^{\bar b_f}$;
\;
$\bar d_i = \bar d^{\bar e^i}$\!\!,\;\;
$\bar d_f 
= \cdots
\bar d_{-1}^{f(-1)}
\bar d_{0}^{f(0)} 
\bar d_{1}^{f(1)}\cdots$,
compare to Section~\ref{SU Defining subgroups by sets of functions}.

\medskip
If for the given $\mathcal X\subseteq \E$ the subgroup $A_{\mathcal X}$ is benign in $F$, then $\bar A_{\mathcal X}
= \langle \bar a_f \;|\; f\in \mathcal X\rangle$ clearly is benign in $\bar F$. If the overgroup $K_{\mathcal X}=\langle
\, X
\;|\;
R
\rangle$ and its subgroup $L_{\mathcal X}\le K_{\mathcal X}$ can explicitly be constructed for $A_{\mathcal X}$, the respective overgroup $\bar K_{\mathcal X}=\langle
\,\bar X
\;|\;
\bar R
\rangle$ and its subgroup $\bar L_{\mathcal X}\le \bar K_{\mathcal X}$ can explicitly be built for $\bar A_{\mathcal X}$ also.

\section{The proof for Theorem~\ref{TH Theorem for rho}}
\label{SU The proof for rho} 

\noindent
For convenience we recall the notation and assumptions of Theorem~\ref{TH Theorem for rho}.

\begin{itemize}
\item[$\circ$]
Let $\mathcal X$ be a set of functions from $\Z$ to $\Z$ with finite support. 

\item[$\circ$]
Let $A_{\mathcal X}$ be the corresponding subgroup of free group $F = \langle a,b,c\rangle$ of rank $3$, see Section~\ref{SU Defining subgroups by sets of functions}. 

\item[$\circ$]
Suppose that $A_{\mathcal X}$ is benign in $F$, i.e., there exist a finitely presented overgroup $K_{\mathcal X}$ of $F$ and a finitely generated subgroup $L_{\mathcal X}$ of $K_{\mathcal X}$ such that $A_{\mathcal X}=K_{\mathcal X} \cap L_{\mathcal X}$. 

\item[$\circ$]
Suppose the groups $K_{\mathcal X}$ and $L_{\mathcal X}$ are given \textit{explicitly}, i.e.:

\begin{itemize}
\item[--]
the group $K_{\mathcal X}$ is given by a finite presentation $K_{\mathcal X}=\langle
\, X
\;|\;
R
\rangle$;

\item[--]
the group $L_{\mathcal X}$ is given by a finite set of generators written as words in the alphabet $X$. 
\end{itemize}
\end{itemize}

For the set $\mathcal Y = \rho \mathcal X$, see Section~\ref{SU The Higman operations}, we want to show that $A_{\mathcal Y}$  is also benign in $F$, and that the corresponding groups $K_{\mathcal Y}$ and $L_{\mathcal Y}$ can be given \textit{explicitly}.
See also Figure~\ref{FI Passage from the benign subgroup} illustrating the passage in Theorem~\ref{TH Theorem for rho}.

\smallskip
The following agreement is going to be useful later:

\begin{Remark}
\label{RE abc can be added}
If the above $K_{\mathcal X}=\langle
\, X
\;|\;
R
\rangle$ and the embedding 
of $F=\langle a,b,c \rangle$ into $K_{\mathcal X}$ are \textit{explicitly} known, it is possible to write the free generators $a,b,c$
as some words $a=a(X),\, b=b(X),\,  c=c(X)$ in the alphabet $X$. In many cases, such as the proofs in \cite{On explicit embeddings of Q}, $K_{\mathcal X}$ is already constructed so that $X$ contains the letters $a,b,c$. But even if $K_{\mathcal X}$ is explicitly given by some other generators \textit{not} involving $a,b,c$, we can apply a Tietze transformation:
add the words $a=a(X),\, b=b(X),\,  c=c(X)$ to the defining relations $R$, and add the letters $a,b,c$ to the generators $X$. Hence we can always assume that the generators of $F$ are among those in $X$. 
\end{Remark}

\subsection{Construction of the direct product $K_P=\bar {\mathcal K} \times \mathscr{A}$}
\label{SU Construction of the direct product K x A}  

Recall that $\mathscr{A}$ is an overgroup of $F$ generated by $a,b,c$, together with the letters $t_1, t'_1, u_1, u_2, d,e$, see \eqref{EQ relations A}. Hence, without loss of generality, we may assume that the latter six letters were not involved in the construction of $K_{\mathcal X}$. 
Since $K_{\mathcal X}$, by construction, contains $a,b,c$, we can assume that $\mathscr{A} \cap K_{\mathcal X}=F$. 
Then the following amalgamated product is correctly defined:
$$
\mathcal K = K_{\mathcal X} *_F \mathscr{A}.
$$
Since $L_{\mathcal X}$ is a subgroup of $K_{\mathcal X}$ such that $F \cap L_{\mathcal X} = A_{\mathcal X}$, we have $\mathscr{A} \cap L_{\mathcal X}   = A_{\mathcal X}$.
That is, $A_{\mathcal X}$ is also benign in the \textit{larger} group $\mathscr{A}$ for the finitely presented $\mathcal K$ and for the same finitely generated subgroup $L_{\mathcal X}$.      

In Section~\ref{NOR SU Auxiliary copy of A built here} we built a copy $\bar{\mathscr{A}}$ of $\mathscr{A}$ on some new generators \eqref{EQ generators of X bar A}. Modifying the steps above we get the copies 
$\bar F,\,
\bar A_{\mathcal X},
\bar {\mathcal K}, 
\bar L_{\mathcal X}$ of the groups $F,\,
A_{\mathcal X},
{\mathcal K}, 
L_{\mathcal X}$, so that $\bar A_{\mathcal X}$ is benign in $\bar{\mathscr{A}}$ for 
$\bar {\mathcal K}$
and 
$\bar L_{\mathcal X}$.

\medskip 

By Remark~\ref{RE finite generated is benign}, 
$\langle a, d, e \rangle$ is benign in $\mathscr{A}$
for the finitely presented $\mathscr{A}$ and for the finitely generated $\langle a, d, e \rangle$.
Hence, the direct product 
$$
P=\bar A_{\mathcal X} \times \langle a, d, e \rangle
$$ is benign in $\bar {\mathscr{A}} \times \mathscr{A}$ for the finitely presented 
$K_P=\bar {\mathcal K} \times \mathscr{A}$ and for the finitely generated $L_P=\bar L_{\!\mathcal X} \!\times \langle a, d, e \rangle \le K_P$.
\medskip 

For each $f\!\in\! \mathcal X$ we by Lemma~\ref{LE action of d_m on f} evidently have
\begin{equation} 
\label{EQ a^d_f = a^d_f in general}
a_f\! = a^{b_f}\! = a^{d_f}\!,
\end{equation}
this simple fact can be explained for, say, $f=(2,5,3)$:
\begin{equation} 
\label{EQ this simple fact can be explained}
a^{d_f}\!
=a^{d_0^{2}\,d_1^5\,d_2^{3}}\!
=(a^{d_0^{2}})^{\,d_1^5\,d_2^{3}}\!
=(a^{b_0^{2}})^{\,d_1^5\,d_2^{3}}\!
=(a^{b_0^{2}\,b_1^5})^{\,d_2^{3}}\!
=a^{b_0^{2}\,b_1^5\,b_2^{3}}\!
=a^{b_f}\! =a_f.
\end{equation}
To verify the steps of \eqref{EQ this simple fact can be explained}, one can first notice that $a^{d_0}$ is noting but the element $a_{f'}^{d_0}$ of Lemma~\ref{LE action of d_m on f} written for the zero sequence $f=(0)$ and for the integer $j=0$. Hence, by that lemma we have $a^{d_0}=a_{f''}=a^{b_0}$ for the sequence $f''=(1)=f'^{+} $.
Next, using Lemma~\ref{LE action of d_m on f} for \eqref{EQ this simple fact can be explained} once again, we have 
$a^{d^2_0}=\big(a^{d_0}\big)^{d_0}=a_{f''}^{d_0}=a_{f'''}$ for the next sequence  $f'''=(2)=f''^{+}$. Applying Lemma~\ref{LE action of d_m on f} several times, we will eventually arrive to the last step of \eqref{EQ this simple fact can be explained}.

Hence,  $A_{\mathcal X} \subseteq \langle a, d, e \rangle$, and, similarly:  
$$
\bar A_{\mathcal X} = \langle \bar a^{\bar b_f} \;|\;  f\in \mathcal X\rangle
= \langle \bar a^{\bar d_f} \;|\;  f\in \mathcal X\rangle\subseteq \langle \bar a, \bar d, \bar e \rangle.
$$
Therefore, the direct product $P$ defined above also lies in $\langle \bar a, \bar d, \bar e \rangle \times \langle a, d, e \rangle$.

\subsection{Obtaining the benign subgroup $Q$}
\label{SU Obtaining the benign subgroup Q}

For each $f\!\in \E$ define in ${\mathscr{A}}$ the auxiliary elements:
\begin{equation*}
\begin{split}
d_{\rho f} = \cdots
d_{1}^{f(-1)}
d_{0}^{f(0)} 
d_{-1}^{f(1)}\cdots 
\quad\;
{\rm and}
\quad\quad
\tilde d_{\rho f} = \cdots
d_{-1}^{f(1)}
d_{0}^{f(0)} 
d_{1}^{f(-1)}\cdots
\end{split}
\end{equation*}
where $\tilde d_{\rho f}$ differs from $d_{\rho f}$ by \textit{reverse order} of its factors $d_{i}^{f(i)}$ only. Say, 
for $f=(2,5,3)$ we have 
$d_{\rho f}
=  d_{-2}^{3} d_{-1}^{5} d_{0}^{2}$
and 
$\tilde d_{\rho f}
=  d_{0}^{2} d_{-1}^{5} d_{-2}^{3}$, compare these with the element $\bar d_{f}
=  \bar d_{0}^{2} \bar d_{1}^{5} \bar d_{2}^{3}$.

In the direct product 
$\bar{\mathscr{A}} \times {\mathscr{A}}$ define the pairs $\lambda_f=\big(\bar d_f,\; \tilde d_{\rho f} \! \big)$. 
The $3$-generator subgroup 
$$
T=\big\langle 
(\bar a,\; a),\;
(\bar d,\; d),\;
(\bar e,\; e^{-1})
\big\rangle
$$ 
clearly contains such 
$\lambda_f$ for every $f\!\in \mathcal E$.
This simple fact requires routine calculations, which are easier to explain via a simple example for, say, $f=(2,5,3)$.
Trivially, $T$ contains the product
$$
(\bar e,\; e^{-1})^{-2} 
\cdot \,
(\bar d,\; d) 
\cdot \,
(\bar e,\; e^{-1})^2
=
\big(\bar d^{\,\bar e^{\;\,2}}
\!\!,\;
d^{e^{\,-2}}\big)
=
\big(\bar d_2
,\;
d_{\,-2}\big),
$$
together with its cube
$
\big(\bar d_2
,\;
d_{\,-2}\big)^3
=
\big(\bar d_2^3
,\;
d_{\,-2}^3\big)
$. 
For similar arguments the product
$$
\lambda_f
=
\big(\bar d_0^2
,\;
d_{0}^2\big)
\cdot
\big(\bar d_1^5
,\;
d_{-1}^5\big)
\cdot
\big(\bar d_2^3
,\;
d_{-2}^3\big)
=
\big( \bar d_{0}^{2} \bar d_{1}^{5} \bar d_{2}^{3},\;\; d_{0}^{2}d_{-1}^{5} d_{-2}^{3} \big)
=\big(\bar d_f,\; \tilde d_{\rho f} \! \big)
$$
is also inside $T$. Hence, the conjugates $(\bar a,\; a)^{\lambda_f}=
\big(\bar a^{\bar d_f},\; a^{\tilde d_{\rho f}}\!\big)$ are also inside $T$  for all $f$.

According to Remark~\ref{RE finite generated is benign}, the $3$-generator group $T$
is benign in $\bar{\mathscr{A}} \!\times\! {\mathscr{A}}$, and for it we can choose the overgroup $K_{T} = \bar{\mathscr{A}} \!\times\! {\mathscr{A}}$ with $L_{T} = T$.
Recall that $P$ is also benign in $\bar{\mathscr{A}} \!\times\! {\mathscr{A}}$. Then, by Corollary~\ref{CO intersection and join are benign multi-dimensional}\;\eqref{PO 1 CO intersection and join are benign multi-dimensional}, the intersection:  
$$
Q = T \cap P
$$
is benign in $\bar{\mathscr{A}} \!\times\! {\mathscr{A}}$ for the finitely presented overgroup:
\begin{equation}
\label{EQ K_Q defined}
\begin{split}
K_Q &
= \big(K_T *_{L_T} v_1\big) 
*_{\bar{\mathscr{A}} \!\times\! {\mathscr{A}}}
\big(K_P *_{L_P} v_2\big)\\
&= \Big(
\big(\bar{\mathscr{A}} \!\times\! {\mathscr{A}}\big) 
*_{T} 
v_1 
\Big)
\, *_{\bar{\mathscr{A}} \times {\mathscr{A}}}
\Big(
\big(\bar {\mathcal K} \!\times\! \mathscr{A}\big)
*_{\bar L_{\!\mathcal X} \times\, \langle a, d, e \rangle} 
v_2 
\Big)
\end{split}
\end{equation}
with two new stable letters $v_1, v_2$, and for the $6$-generator subgroup $L_Q\!=\big(\!\bar{\mathscr{A}} \!\times\! {\mathscr{A}}\big)^{v_1v_2}$\!.

\medskip
Our next aim is to describe a structure of $Q$. First, we note that $Q$ lies in the direct product $P=\bar A_{\mathcal X} \times \langle a, d, e \rangle$. 
Thus, each element $q \in Q$ can be written in the form $q = (q_1,q_2)$, where 
$q_1 \in \bar A_{\mathcal X}$ and
$q_2 \in \langle a, d, e \rangle$.
The first coordinate $q_1$ is a word in some elements $\bar a_{f_i}$, say, $q_1 = w(\bar a_{f_1},\ldots,\bar a_{f_n})$. 
By analogy with \eqref{EQ a^d_f = a^d_f in general}, we have:
$$
\bar a_f\! = \bar a^{\bar b_f}\! = \bar a^{\bar d_f}.
$$
We can apply this for each of the above $f_i$ and, therefore, $q_1$ can be rewritten as a word in three letters $\bar a, \bar d, \bar e$, say, $q_1 = v(\bar a, \bar d, \bar e)$.

On the other hand, since $Q$ lies in $T=\big\langle 
(\bar a,\; a),\;
(\bar d,\; d),\;
(\bar e,\; e^{-1}) 
\big\rangle$, 
we have $q_2 = v(a, d, e^{-1})$.
From this and from $q_1 = w(\bar a_{f_1},\ldots,\bar a_{f_n})
= w(\bar a^{\bar d_{f_1}},\ldots,\bar a^{\bar d_{f_n}})$ we easily conclude that 
$q_2 = w(a^{\tilde d_{\rho f_1}},\ldots, a^{\tilde d_{\rho f_n}})$.

Let us explain this with the example of $q$ such that $q_1 = \bar a_f$ for $f = (2,5,3)$. 
Then the first coordinate is:
\begin{equation*}
\begin{split}
q_1 &
= \bar a_f 
= \bar a^{\bar b_f}
=\bar a^{\bar d_f}
= \bar a^{\bar d_{0}^{2} \bar d_{1}^{5} \bar d_{2}^{3}}
= \bar d_{2}^{\;-3}\bar d_{1}^{\;-5}  \bar d_{0}^{\;-2}\, \cdot \,
\bar a \,
\cdot\, \bar d_{0}^{2}\, \bar d_{1}^{5}\, \bar d_{2}^{3}\\
&= 
\bar e^{\;-2} \,\bar d^{\,-3} \,\bar e^{\;2} 
\cdot 
\bar e^{\;-1} \,\bar d^{\,-5} \,\bar e
\cdot 
\bar d^{-2}
\,\cdot\;\, 
\bar a 
\,\; \cdot \; \bar d^{2}
\cdot
\bar e^{\;-1} \,\bar d^{\,5} \,\bar e
\cdot 
\bar e^{\;-2} \,\bar d^{\,3} \,\bar e^{\;2}
\end{split}
\end{equation*}
with respect to which the second coordinate has to be: 
\begin{equation*}
\begin{split}
q_2& =
e^{\;2} \, d^{\,-3} \, e^{\;-2} 
\cdot \,
e \, d^{\,-5} \, e^{\;-1}
\cdot \,
d^{-2}
\cdot\,\,
a \,\,
\cdot \, d^{2}
\cdot\,
e \, d^{\,5} \, e^{\;-1}
\cdot \,
e \, d^{\,3} \, e^{\;-2}\\
& = 
d_{-2}^{\;-3} d_{-1}^{\;-5}   d_{0}^{\;-2}\, \cdot\,  \,
a \,
\, \cdot\,  d_{0}^{2}\,  d_{-1}^{5}\,  d_{-2}^{3}
= a^{d_{0}^{2}\,  d_{-1}^{5}\,  d_{-2}^{3}}=a^{\tilde d_{\rho f}}.
\end{split}
\end{equation*}
Thus the benign subgroup $Q$ of $\bar{\mathscr{A}} \!\times\! {\mathscr{A}}$, in fact, has simple description:
$$
Q = \big\langle 
\big(\bar a^{\bar d_f},\; a^{\tilde d_{\rho f}}\!\big)
\;|\;  f\in \mathcal X
\big\rangle
= \big\langle 
(\bar a,\; a)^{\lambda_f}
\;|\;  f\in \mathcal X
\big\rangle.
$$
This in particular means that $Q$ lies inside $\bar F \!\times\! F$, and so $Q$ is benign in $\bar F \!\times\! F$ also for \textit{the same} choice of $K_Q$ and $L_Q$ made earlier.
%

\subsection{``Extracting'' $A_{\rho \mathcal X}$ from $Q$}
\label{SU Extracting A rho X from Q} 

Our next objective is to modify the obtained benign subgroup $Q$ via a few steps  ``to extract'' the benign subgroup 
$A_{\rho \mathcal X}=
\big\langle a^{b_{\rho f}}
\;|\;  f\in \mathcal X
\big\rangle
$ from it.

Using Remark~\ref{RE order of d_i does not matter} and equation \eqref{EQ a^d_f = a^d_f in general}, we deduce
$$
a^{\tilde d_{\rho f}}
=a^{d_{\rho f}}
=a^{b_{\rho f}}
=a_{\rho f},
$$
that is, $A_{\rho \mathcal X}$ is nothing but the group generated by the second coordinates $a^{\tilde d_{\rho f}}$ of the above mentioned pairs $q = (q_1,q_2)$ from $Q$.

$\bar F \cong \bar F \times \1$ is benign in 
$\bar F \times F$ for the finitely presented $\bar F \times F$ and finitely generated $\bar F \times \1$, see Remark~\ref{RE finite generated is benign}.  
Hence, by Corollary~\ref{CO intersection and join are benign multi-dimensional}\;\eqref{PO 2 CO intersection and join are benign multi-dimensional} the join 
$$
Q_1=\big\langle \bar F \!\times\! \1 ,\, Q \big\rangle
=\bar F \!\times\! \langle
a_{\rho f} \mathrel{|} 
f \in \mathcal X \rangle
$$ 
is benign in $\bar F \!\times\! F$ for the finitely presented overgroup:
$$
K_{Q_1}= 
\big((\bar F \!\times\! F)\,*_{\bar F \times \1} w_1\big) 
\,*_{\bar F \times F}
\big(K_{Q} *_{L_{Q}} w_2\big)
$$
with two further new letters $w_1, w_2$,\; and for its $12$-generator subgroup:
$$L_{Q_1} = \big\langle(\bar F \!\times\! F)^{w_1}\!,\;(\bar F \!\times\! F)^{w_2} \big\rangle.$$

Finally,  
$F = \1 \!\times\!  F $ is benign in 
$\bar F \!\times\! F$ for the finitely presented $\bar F \!\times\! F$ and finitely generated $\1 \!\times\!  F $.  
Hence, 
by Corollary~\ref{CO intersection and join are benign multi-dimensional}\;\eqref{PO 1 CO intersection and join are benign multi-dimensional}
the intersection:
$$
A_{\rho \mathcal X}=A_{\mathcal Y}=
\big(\1 \!\times\!  F \big) \cap Q_1 
=
\langle
a_{\rho f} \mathrel{|} 
f \in \mathcal X \rangle
$$ 
is benign in $\bar F \!\times\! F$ for the finitely presented:
$$
K_{\rho \mathcal X}
=
K_{\mathcal Y}
=
\big((\bar F \!\times\! F)\,*_{\1 \times  F} w_3\big) 
\,*_{\bar F \times F}
(K_{Q_1} *_{L_{Q_1}} w_4)
$$
with stable letters $w_3, w_4$,\, and for the $6$-generator subgroup
$$
L_{\rho \mathcal X}
=L_{\mathcal Y}
=(\bar F \!\times\! F)^{w_3 w_4}.
$$
But since $A_{\rho \mathcal X}$ is inside $F$, it is benign in $F$ also, for the \textit{same} choice of $K_{\rho \mathcal X}$ and $L_{\rho \mathcal X}$ above.

\subsection{Writing $K_{\rho \mathcal X}$ by its generators and defining relations explicitly}
\label{SU Writing K rho X  by its generators and defining relations}

From sections \ref{SU The group A} and \ref{NOR SU Auxiliary copy of A built here} we know some finite presentations of the groups $\mathscr{A}$ and  $\bar {\mathscr{A}}$, see notation in \eqref{EQ generators of XA} and in \eqref{EQ generators of X bar A}:
$$
\text{$\mathscr{A}=\langle
\,X_{\! \mathscr{A}}
\;|\;
R_{\! \mathscr{A}}
\rangle$
\;and\;
$\bar {\mathscr{A}}=\langle
\,X_{\! \bar {\mathscr{A}}}
\;|\;
R_{\! \bar {\mathscr{A}}}
\rangle$.}
$$
Also recall that the groups 
$K_{\mathcal X}=\langle
\,X
\;|\;
R
\rangle$ and $\bar K_{{\mathcal X}}=\langle
\,\bar X
\;|\;
\bar R
\rangle$ together with the finitely many generators of $L_{\!\mathcal X}$ and of $\bar L_{\!\mathcal X}$
are supposed to be explicitly given.
By Remark~\ref{RE abc can be added} we may assume  
$X_{\! \mathscr{A}} \cap\, X = \{a,b,c\}$ and 
$X_{\! \mathscr{\bar A}} \cap\, \bar X = \{\bar a,\bar b,\bar c\}$.

The amalgamated product $\mathcal K = K_{\mathcal X} *_F \mathscr{A}$
can be generated by the generators 
$X \backslash \{a,b,c\}$ of $K_{\mathcal X}$, together with the generators $X_{\! {\mathscr{A}}}$ (we exclude $a,b,c$ from $X$ because they already were included in $X_{\! {\mathscr{A}}}$, see Remark~\ref{RE abc can be added}).
As defining relations of $\mathcal K$ we may take $R \cup R_{\! \mathscr{A}}$. We can similarly treat the copy $\bar {\mathcal K}$ of ${\mathcal K}$.

Then the direct product $K_P=\bar {\mathcal K} \times \mathscr{A}$
can be given by the above mentioned generators and relation, \textit{plus} the relations making sure $\bar {\mathcal K}$ commutes with $\mathscr{A}$:
\begin{equation} 
\label{EQ K_P presentation}
\begin{split}
K_P=
\Big\langle
X_{\! \mathscr{A}},\;  X_{\! \bar{\mathscr{A}}}\,;\;\;
\bar X 
\,\backslash\, \{ \bar a, \bar b, \bar c \}\;
\;\mathrel{|} \;\; 
R_{\mathscr{\!A}};\;
R_{\mathscr{\!\bar A}};\;
\bar R;\\
&\hskip-59mm 
\text{generators in $X_{\! \mathscr{A}}$ commute with those in $X_{\! \bar{\mathscr{A}}}$}\,;\; \\
&\hskip-59mm 
\text{generators in $X_{\! \mathscr{A}}$ commute with those in $\bar X 
\,\backslash\, \{ \bar a, \bar b, \bar c \}$}
\Big\rangle\,.
\end{split}
\end{equation}

Next, taking into account the way we constructed $K_{T}$, $K_{Q}$, $K_{Q_1}$, $K_{\rho \mathcal X}$ with \textit{fixing} effect of the new letters $v_1, v_2;\; w_1, w_2, w_3, w_4$ (on respective finitely generated subgroups) we have
\begin{equation}
\label{EQ K rho X}
\begin{split}
K_{\rho \mathcal X}=
\Big\langle
X_{\! \mathscr{A}},\;  X_{\! \bar{\mathscr{A}}};\;\;
\bar X_{\! \mathcal X} 
\backslash \{ \bar a, \bar b, \bar c \};\;
v_1, v_2;\;  w_1, w_2, w_3, w_4 
\;\mathrel{|} \;\; 
R_{\mathscr{A}};\;
R_{\mathscr{\bar A}};\;
\bar R;\\
&\hskip-99mm 
\text{generators in $X_{\! \mathscr{A}}$ commute with those in $X_{\! \bar{\mathscr{A}}}$};\; \\
&\hskip-99mm 
\text{generators in $X_{\! \mathscr{A}}$ commute with those in $\bar X 
\,\backslash\, \{ \bar a, \bar b, \bar c \}$};\; \\
&\hskip-99mm 
\text{$v_1$ fixes $\bar a a,\;
\bar d d,\;
\bar e e^{-1} $};
\;\;\;
\text{$v_2$ fixes $\bar L_{\!\mathcal X}$ and $a,d,e$};\; 
\\
&\hskip-99mm \text{$w_1$ fixes  
$\bar a,\bar b,\bar c$};
\;\;
\text{$w_2$ fixes  $X_{\!\! \mathscr{A}}^{v_1 v_2}$ and $X_{\!\! \bar{\mathscr{A}}}^{v_1 v_2}$};\;\;
\text{$w_3$ fixes $ a, b, c$}
\\
&\hskip-99mm  
\text{$w_4$ fixes 
$\big\{
a,b,c, \bar a, \bar b, \bar c
\big\}^{\! w_1} \! \cup \big\{
a,b,c, \bar a, \bar b, \bar c
\big\}^{\! w_2}$}
\Big\rangle\,.
\end{split}
\end{equation}
As $L_{\rho \mathcal X}$ we can explicitly take the $6$-generator subgroup
$\big\langle
a,b,c, \bar a, \bar b, \bar c
\big\rangle^{w_3 w_4}$\! in $K_{\rho \mathcal X}$ by Corollary~\ref{CO intersection and join are benign multi-dimensional}\;\eqref{PO 1 CO intersection and join are benign multi-dimensional}.
%
%
In \eqref{EQ K rho X}\; ``$w_2$ fixes  $X_{\!\! \mathscr{A}}^{v_1 v_2}$ and $X_{\!\! \bar{\mathscr{A}}}^{v_1 v_2}$'' means that conjugation by $w_2$ fixes 
the conjugates of \textit{each} of the generators $X_{\!\! \mathscr{A}}$ and $X_{\!\! \bar{\mathscr{A}}}$
from \eqref{EQ generators of XA} and \eqref{EQ generators of X bar A} by the product $v_1 v_2$.

If $K_{\mathcal X}$ has $m$ generators (which we may assume include $a,b,c$) and $n$ defining relations, and if 
$L_{\mathcal X}$ has $k$ generators, then the group $K_{\rho \mathcal X}$ in \eqref{EQ K rho X} has 
$9 + 9 + (m - 3) + 2 +4 = m +21$ generators and 
$20+20+n+9\cdot 9 + 9\cdot (m-3)
+3 + k +3 + 3  + 2\cdot 9 +3 +2\cdot 6 
= n+9m+k+136$ 
defining relations.

\medskip
In the above constructions we have supposed that the overgroup $K_{\mathcal X}$ of $F$ has $a,b,c$ among its generators by Remark~\ref{RE abc can be added}. Also observe a formatting issue in \eqref{EQ K rho X}: we write not 
``$v_1$ fixes
$(\bar a,\; a)$,
$(\bar d,\; d)$,
$(\bar e,\; e^{-1})$''
but 
``$v_1$ fixes
$\bar a a$,
$\bar d d$,
$\bar e e^{-1}$''
which has the same meaning  because all generators in 
$X_{\! \mathscr{A}}$ already commute with those in $X_{\! \bar{\mathscr{A}}}$ according to the second line of \eqref{EQ K rho X}.

\medskip
Let us finish this chapter with an example promised in \textit{\nameref{SE Introduction}}, to show why Higman's theorem is not an \textit{if and only if} theorem for \textit{all} countably generated groups, unlike the case of finitely generated groups.

\begin{Example}
\label{EX Why Higman's theorem is not iff for countably generated groups}
For any subset $\Pi$ of the set of all prime numbers denote by $\mathbb Q_\Pi$ the additive group of all rationals with denominators divisible by primes from the subset $\Pi$ only. As the recursive group $\mathbb Q$ is embeddable into a finitely presented group, all subgroups $\mathbb Q_\Pi$ of $\Q$ have to be  embeddable, as well. There are uncountably many such subsets $\Pi$, and hence, such groups $\mathbb Q_\Pi$ also are uncountably many. But the set of all recursive groups is at most countable, and so ``many'' of the groups $\mathbb Q_\Pi$ cannot be recursive. 
\end{Example}

\section{Appendix: The main steps of Higman's construction}
\label{SU Some of the main steps of Higman's construction} 

\noindent
In this appendix we would like to outline some of the main steps of Higman's embedding in \cite{Higman Subgroups of fP groups} to identify its fragments to which the current article concerns, and which we attempt to make explicit.

Notice that in the sections above we used $F$ to denote both the free groups of \textit{any} finite rank, as well as, the specific free group $F=\langle a,b,c\rangle$ of rank $3$, like in Theorem~\ref{TH Theorem for rho}. In the current appendix within the same sentences we are going to use
$F$ to denote the free group of \textit{arbitrary} finite rank, its particular instance $F_3=\langle a,b,c\rangle$ of rank $3$, and also $F_2=\langle x,y\rangle$ of rank $2$. Hopefully, this distinction in notation will make it clear which of the free groups $F,\; F_3, F_2$ we mean.

\subsection{Translating recursion to the language of benign subgroups in free groups of arbitrary finite rank}
\label{SU Translating recursion to the language} 

Higman begins by the classic Kleene characterization of recursively enumerable subsets of non-negative integers $\N_0$. 
Then the G{\" o}del numbering allows to generalize the defined notion of recursively enumerable subsets $\mathcal X$ of $\N_0$ to the subsets $\mathcal X$ of \textit{any} effectively enumerable set $\E$, in particular, of the set $\E$ of all functions $f\! : \Z \to \Z$ with finite supports, that we mentioned in Section~\ref{SU Integer functions f}. 

Theorem 3 in \cite{Higman Subgroups of fP groups} states that the recursively enumerable subsets $\mathcal X$ of the above functions set $\E$ are exactly those that 
can be constructed from two specific  subsets 
$\Zz$ and
$\S$ by means of the Higman ope\-ra\-tions \eqref{EQ Higman operations}, see sections \ref{SU Integer functions f} and \ref{SU The Higman operations} above. 
The set of all such subsets of $\E$ is denoted by $\mathscr{S}$\!.

In the free group $F_3=\langle a,b,c\rangle$, using the sets $\mathcal X$, the respective specific subgroups  $A_{\!\mathcal X} =\langle
a_f \;|\; f\in \mathcal X
\rangle$ are defined.
Theorem~4 in \cite{Higman Subgroups of fP groups} states that $\mathcal X$ is recursively enumerable in $\E$ if and only if the respective subgroup $A_{\!\mathcal X}$ is \textit{benign} in $F_3$, see sections \ref{SU Defining subgroups by sets of functions} and \ref{SU Benign subgroups} above. This so far translates recursion to the language of benign subgroups of very \textit{specific type} in the free group of rank $3$. 

Lemma 5.1 in \cite{Higman Subgroups of fP groups} establishes connection between \textit{all} recursively enumerable subsets in the free group of rank $2$, and benign subgroups of a certain \textit{specific type} in the free group of rank $3$.
Then Lemma 5.2 in \cite{Higman Subgroups of fP groups} shows that in a free group of \textit{arbitrary} finite rank, a subgroup is recursively enumerable if and only if it is benign.

\medskip
The above steps translate the logical notion of \textit{recursion} to the group theoretic concept of \textit{benign subgroups} of free groups of \textit{any} finite rank, and the rest of discussion can be conducted in the language of benign subgroups.

\subsection{Outputting the recursive set $\mathcal X$}
\label{SU Outputting the recursive set X} 
Only after the above outlined translation from the recursively enumerable subsets to the language of benign subgroups, \cite{Higman Subgroups of fP groups} turns to the recursive group $G$ in question, assuming it has the presentation $G = \langle\, X \mathrel{|} R \,\rangle$ with a finite set of generators $X$ and with a recursively enumerable set of defining relations $R$, that is, we have $G \cong F/\langle R \rangle^F$  for some free group $F$ of rank $|\,X|$. 

As the set $R$ is recursively enumerable, its normal closure $\langle R \rangle^F$ also is recursively enumerable for very simple combinatorial reasons (in \cite{Higman Subgroups of fP groups} this closure is denoted by $R$).

The group $F$ can be embedded into a free group $F_2=\langle x,y \rangle$ of rank $2$, see \textit{p.}\,474 in \cite{Higman Subgroups of fP groups}, using any of the well known textbook methods, see \cite{Kargapolov Merzljakov,Bogopolski,Robinson}. 
By the remark preceding Lemma~3.8 in \cite{Higman Subgroups of fP groups}, the normal subgroup $\langle R \rangle^F$ is benign in $F$ if and only if it is benign in $F_2$. 
Hence, we can imagine all the relations from $\langle R \rangle^F$ are rewritten in just \textit{two} letters $x,y$.
The purpose of this passage is that it allows to trivially output the integer sequences $f\in \mathcal X$ from words on $x,y$. For example, from the relation $x^2 y^5 x^3 \in F_2$ one can extract the sequence $f\!=\!(2, 5, 3)\in \mathcal X$, which later is being used in the free group $F_3=\langle a,b,c\rangle$ to produce the element 
$b_f=
b_{0}^{2} 
b_{1}^{5}
b_{2}^{3}$, and then the conjugate $a_{f}\!=a^{b_f}$ inside $A_{\mathcal X}$, 
see Section~\ref{SU Defining subgroups by sets of functions} above.

In the next steps of \cite{Higman Subgroups of fP groups} Higman is going to use this benign subgroup $A_{\mathcal X}$ of $F_3$ 
to build, for the benign subgroup $\langle R \rangle^F$ of $F_2$, the respective finitely presented overgroup $K_{\mathcal X}$ with a finitely generated subgroup $L_{\mathcal X} \le K_{\mathcal X}$ such that $F_2 \cap L_{\mathcal X} = \langle R \rangle^F$\! holds.
These steps, unfortunately, create two hurdles for the future.

\smallskip
\textit{Firstly}, Higman's process outputs the groups $K_{\!\mathcal X}$ and  $L_{\!\mathcal X}$ not for the initial free group $F$ but for the auxiliary group $F_2$ used above.  Thus, following \cite{Higman Subgroups of fP groups}, we have to ``return back'' from $F_2$ to $F$ in order to modify the obtained $K_{\!\mathcal X}$ and  $L_{\!\mathcal X}$ for the free group $F$ of rank $|\,X|$ from the original presentation $G=F/\langle R \rangle^F$\!.
That part is not hard to do by a few routine steps because by the remark before Lemma~3.8 in \cite{Higman Subgroups of fP groups}, a group is benign in $F_2$ if and only if it is benign in $F$.
In order \textit{not} to introduce extra new notation, suppose the above $K_{\mathcal X}$ and $L_{\mathcal X}$ are those found for $F$, that is, $F \cap L_{\mathcal X} = \langle R \rangle^F$ holds inside a certain finitely presented $K_{\mathcal X}$ for the recursive subgroup $\langle R \rangle^F$\!.

\smallskip
\textit{Secondly}, the passage from $R$ to $\langle R \rangle^F$ very much expands the set of relations we have to work with. For, even if the group $G = \langle\, X \mathrel{|} R \,\rangle$ is defined with some very simple set $R$, the subgroup $\langle R \rangle$ and the normal closure $\langle R \rangle^F$ may contain a by far wider stockpile of relations, which generate a by far wider set of functions $\mathcal X$, and make it by far harder to write it via the operations 
\eqref{EQ Higman operations} 
in the coming steps.

\begin{Remark}
\label{RE we do not use <R>^F}
Notice that in \cite{On explicit embeddings of Q} we were able to avoid the routine of these two hurdles by using other methods that extract a by far smaller  set $\mathcal X$ by means of a ``preliminary'' overgroup $T_\Q$ of $\Q$. 
\end{Remark}

Anyway, in order to be able to proceed to Higman's next step, assume the set $\mathcal X$ is given somehow, be it some description via Kleene characterization, or via a Turing machine, etc.

\subsection{Building $\mathcal X$ via the operations \eqref{EQ Higman operations}}
\label{SU Buildinng X via H} 

The next objective is to \textit{explicitly construct} the known set $\mathcal X$ via the Higman operations
\eqref{EQ Higman operations} from $\Zz$ and $\S$.
Lemma 2.8 in \cite{Higman Subgroups of fP groups} provides a long but yet \textit{explicit} algorithm to construct certain partial recursive functions $f(n,r)$, $a(r)$, $b(r)$, with $r=0,1,2,\ldots$, which one-by-one write down all the sequences $g\in\mathcal X$, and indicate their start- and end points, in the sense that for any such $g$ there is some $r$ for which $g(n)=f(n,r)$ for all $n\in \Z$, while $f(n,r)=0$ for all $n<a(r)$ and $n>b(r)$. 
Such functions 
$f_{\E}(n,r)$, 
$a_{\E}(r)$, 
$b_{\E}(r)$  
are explicitly built first for the case of the whole set $\mathcal X = \E$. Then for any other generic recursively enumerable subset $\mathcal X \subseteq \E$ there exists a partial recursive function $h(s)$ such that the $r$'th function $f_{\E}(n,r)$ corresponds to a function in $\mathcal X$ if and only if $r=h(s)$ for some $s$. 

For every multi-variable integer-based function $f(x_1,\ldots,x_n)$ its \textit{graph} is being defined before Lemma 2.2, and $\mathscr F$ is  denoted to be the set of all functions whose graphs are in $\mathscr S$. Lemmas 2.1\,--\,2.7in \cite{Higman Subgroups of fP groups} show that $\mathscr F$ contains all the partial recursive functions. 
Hence, the graphs of the above functions $f,a,b$ also are in $\mathscr F$, i.e., their graphs can be constructed via the operations \eqref{EQ Higman operations}.
Using the above steps, Higman finishes the proof of Theorem 3 on pages 463\,--\,464 by showing how the earlier picked sequence $g\in \mathcal X$ can constructively be obtained  via \eqref{EQ Higman operations} using the graphs of these functions $f,a,b$.  

\smallskip 
Although each of the steps with $f,a,b$ listed above is \textit{explicit}, and is theoretically doable, they require such a vast routine of operations that their application to non-trivial groups is a very hard task. 

\begin{Remark}
\label{RE we do not use f a b}
In our construction of explicit embeddings for $\Q$ in \cite{On explicit embeddings of Q} we never follow this portion of Higman's steps with $f,a,b$. Instead we use the much shorter constructions from \cite{The Higman operations and  embeddings} that let us build $\mathcal X$ from \eqref{EQ Higman operations} for our groups, such as $\Q$.
\end{Remark}

\subsection{Explicit construction of $K_{\!\mathcal Y}$ and  $L_{\!\mathcal Y}$ for all \eqref{EQ Higman operations}}
\label{SU Explicit construction of KX and LX}

Next, assume the set $\mathcal X$ has already been constructed from $\Zz$ and
$\S$ via some \textit{known} sequence of operations \eqref{EQ Higman operations}. 

By Theorem~4 in \cite{Higman Subgroups of fP groups}, the construction of $\mathcal X$ via \eqref{EQ Higman operations} is possible if and only if $A_{\!\mathcal X}$ is benign in $F=F_3=\langle a,b,c \rangle$, that is, the free group $F_3$ is embeddable into a finitely presented group $K_{\!\mathcal X}$ with a finitely generated subgroup  $L_{\!\mathcal X}$ such that $F_3 \cap L_{\!\mathcal X} = A_{\!\mathcal X}$ holds. 

As we mentioned in Section~\ref{SU HMotivation for Higman's reversing operation rho and Theorem 1.1}, the current step will become explicit, if for \textit{each} of the operations \eqref{EQ Higman operations} we are able to explicitly write the groups $K_{\!\mathcal X}$ and $L_{\!\mathcal X}$ for two initial sets 
$\mathcal X = \Zz, \S$, and then we are in position to explicitly write the groups $K_{\!\mathcal Y}$ and $L_{\!\mathcal Y}$ as soon as the set $\mathcal Y$ is obtained from a certain set $\mathcal X$ by one of operations \eqref{EQ Higman operations}, and the groups $K_{\!\mathcal X}$ and $L_{\!\mathcal X}$ already are explicitly known for that $\mathcal X$. 

The purpose of Theorem~\ref{TH Theorem for rho} of the current paper is to accomplish  the stated job for Higman's reversing operation $\rho$ from the list \eqref{EQ Higman operations}.

\subsection{Passage from the benign subgroup $A_{\!\mathcal X}$ in $F_3$ to the benign subgroup $\langle R \rangle^F$ in $F$}
\label{SU Passage from the benign subgroup} 

After two groups $K_{\!\mathcal X}$ and  $L_{\!\mathcal X}$ are found for the benign subgroup $A_{\!\mathcal X}$ of $F_3=\langle a,b,c \rangle$, Higman constructs the analogs of these two groups for the benign subgroup $\langle R \rangle^F$ of $F_2$. 
This is done within constructions of Lemma 5.1 and Lemma 5.2 in \cite{Higman Subgroups of fP groups}.

Further, as we warned in Section~\ref{SU Outputting the recursive set X} above, after $K_{\!\mathcal X}$ and  $L_{\!\mathcal X}$ are found for the benign subgroup $\langle R \rangle^F$ of $F_2$, Higman has to ``return back'' from $F_2$ to the initial free group $F$ from the presentation presentation $G=F/\langle R \rangle^F$\!,\, to find the respective finitely presented overgroup $K_{\mathcal X}$ with a finitely generated subgroup $L_{\mathcal X} \le K_{\mathcal X}$ such that $F \cap L_{\mathcal X} = \langle R \rangle^F$\! holds.

Notice that to avoid complication of the notation, we use the same symbols $K_{\mathcal X}$ and $L_{\mathcal X}$ to denote them. 
As we mentioned in Remark~\ref{RE we do not use <R>^F}, in our constructions in \cite{On explicit embeddings of Q} we do not have to use these passages from $F$ to $F_2$, and then from $F_2$ back to $F$.

Anyway, to proceed to the last step assume $K_{\!\mathcal X}$ and  $L_{\!\mathcal X}$ have been found for the normal benign subgroup $\langle R \rangle^F$ in the \textit{initial} free group $F$.

\subsection{The group $K$ and the ``Higman Rope Trick''}
\label{SU The group K and the Higman Rope Trick} 

For the free group $F$ it is easy to pick such an isomorphic copy $F'$ of it, such that $F'$ intersects with $F$ precisely in $\langle R \rangle^F$\! inside a larger common overgroup. For example, one can first embed $F$ into the HNN-extension $F *_{\langle R \rangle^F} r$, and then choose $F'$ to be the conjugate $F^r$ of $F$ by the stable letter $r$. 

Lemma~3.5 in \cite{Higman Subgroups of fP groups} provides an alternative  definition for benign subgroups. Applying it to  $\langle R \rangle^F$\! one gets that it is benign in $F$ if and only if the free product $H= F *_{\langle R \rangle^F} F'$ of the group $F$ and of its copy $F'$, amalgamated in $\langle R \rangle^F$\!,  is embeddable into some finitely presented group $K$. From the simple proof of Lemma~3.5 it is clear that as $K$ one could choose the HNN-extension $K = K_{\mathcal X} *_{L_{\mathcal X}} s$ of the above $K_{\mathcal X}$ with the fixed  subgroup $L_{\mathcal X}$
(where $K_{\mathcal X}$ and $L_{\mathcal X}$ are those outputted at the end the previous step \ref{SU Passage from the benign subgroup}). This extension is generated by \textit{finitely} many generators of $K_{\mathcal X}$ together with the stable letter $s$. Also, $K$ is \textit{finitely} defined by the finite set of defining relations of $K_{\mathcal X}$ together with the finitely many relations telling that $s$ fixes the finitely many generators of $L_{\mathcal X}$. 

The embedding of $H$ into $K$ can trivially be continued to the embedding of $H$ into the direct product $K \times G$, for which we in addition to finitely many relations of $K$ may have \textit{infinitely} many relations for $G$. Plus, the \textit{finitely} many relations telling that the finitely many generators of $K$ and $G$ commute, of course. 

The free group $F$ has two isomorphic images $F$ and $F'=F^r$ in the \textit{first} direct factor of $K\! \times \! G$, and $F$ has its surjective image $F/\langle R \rangle^F = G$ as the \textit{second} direct factor. 

For any word $f\!\in F$ we inside $K \times G$ can define 
the identity map $(f^r, 1)  \to (f^r, 1)$, and the map
$(f, 1) \to (f,\; \langle R \rangle^F\! f)$, with the coset $\langle R \rangle^F\! f$ from $F/\langle R \rangle^F = G$. Whenever the word $f$ is a \textit{relation} of $G$, we have $\langle R \rangle^F\! f = \langle R \rangle^F\! = 1$ in $G$, that is, these two maps \textit{agree} on the amalgamated subgroup $\langle R \rangle^F$\!, and they can therefore be continued on the whole $H$, to form an injective homomorphism $\alpha$ from $H$ into $K \times G$.

Hence one can define the last HNN-extension $(K \times G) *_{\alpha} t$ which turns out to be the desired finitely presented overgroup of $G$. This is shown by an elegant process often called the ``Higman Rope Trick'' that proves that all but \textit{finitely} many  of the relations of this HNN-extension are redundant. 

Namely, the action of $t$ on $K \times G$ is defined by the images of $\alpha$ on \textit{finitely} many generators of $F$ and of its copy $F'$.
Thus, the only remaining spot with infinitely many defining relations is $G$, and Higman's charming trick is to unveil why they may simply be ignored, 
also see the helpful discussion \cite{Higman rope trick}.

\subsection{Recapping the summary}
\label{SU Recapping the summary} 
To summarize the review in this appendix, we observe that: some parts of Higman's original construction in \cite{Higman Subgroups of fP groups} already are constructive, while others are not. 
Hence, Theorem~\ref{TH Theorem for rho} of the current paper closes one of the remaining gaps that was not covered by an explicit constructions yet. 
This theorem is a small step towards the bigger objective of obtaining the constructive version of the Higman Embedding Theorem.

\medskip
\noindent 
E-mail:
\href{mailto:v.mikaelian@gmail.com}{v.mikaelian@gmail.com}
$\vphantom{b^{b^{b^{b^b}}}}$

\noindent 
Web: 
\href{https://www.researchgate.net/profile/Vahagn-Mikaelian}{researchgate.net/profile/Vahagn-Mikaelian}

\end{document}